\documentclass{amsart}
\usepackage{amssymb,latexsym,amsmath,amsfonts}

\numberwithin{equation}{section}
\theoremstyle{definition}
\newtheorem{definition}{Definition}[section]
\newtheorem{remark}[definition]{Remark}
\theoremstyle{plain}
\newtheorem{theorem}[definition]{Theorem}
\newtheorem{lemma}[definition]{Lemma}
\newtheorem{proposition}[definition]{Proposition}
\newtheorem{example}[definition]{Example}

\newcommand{\Cn}{{\mathbb{C}}^n}

\newcommand{\bdy}{\partial \Omega}
\newcommand{\BD}{\partial \omega}
\newcommand{\mn}{\mathcal{M}}
\newcommand{\MN}{\mathfrak{m}}
\newcommand{\al}{|\alpha|}
\newcommand{\bet}{|\beta|}
\newcommand{\er}{{\mathfrak{r}}{\mathfrak{e}}}
\newcommand{\mi}{{\mathfrak{i}}{\mathfrak{m}}}

\newcommand{\Fie}{\varPhi}
\newcommand{\Rm}{{\mathbb{R}}^m}
\newcommand{\Il}{\mathbb{L}}
\newcommand{\SIl}{\text{\bf S}}
\newcommand{\prr}{(w_*,w_n)}
\newcommand{\vect}[2]{\dfrac{\partial #2}{\partial w_{#1}}}
\newcommand{\Vct}[1]{\dfrac{\partial}{\partial z_{#1}}}

\newcommand{\cvct}[2]{\dfrac{\partial #2}{\partial {\bar w}_{#1}}}
\newcommand{\gamder}[2]{\dfrac{\partial{\gamma}_{#1}}{\partial x_{#2}}}
\newcommand{\fuder}[2]{\dfrac{\partial (\overline{F_{#1}}\circ\gamma)}{\partial x_{#2}}}
\newcommand{\peeq}{{\mathcal{P}}_q}
\newcommand{\rng}{q \in W(p)\cap \text{Image}(\gamma)}
\newcommand{\Al}{\overline{A_{\mu}(w})}
\newcommand{\jA}{A_j(w)}
\newcommand{\Fl}{\overline{F_{\mu}(w})}
\newcommand{\jF}{F_j(w)}
\newcommand{\dAl}[1]{\dfrac{\partial\overline{A_{\mu}}}{\partial w_{#1}}(w)}
\newcommand{\dlF}[1]{\dfrac{\partial\overline{F_{\mu}}}{\partial w_{#1}}(w)}
\newcommand{\zum}{\sum_{\mu,j=1}^{n-1}}
\newcommand{\eps}{\varepsilon}

\newcommand{\lgz}{\log(z{\bar z})}
\newcommand{\levi}[1]{{\mathfrak{L}}^{(#1)}}

\newcommand{\boo}{\text{\bf v}}
\newcommand{\Am}{\Im\mathfrak{m}}
\newcommand{\ws}{\varrho}

\begin{document}
\title[Analytic Interpolation Manifolds]{On Analytic Interpolation Manifolds in \\
Boundaries of Weakly Pseudoconvex Domains}
\author{Gautam Bharali}
\address{Mathematics Department, University of Wisconsin-Madison, 480 Lincoln Drive,
Madison, WI 53706}
\email{bharali@math.wisc.edu}
\keywords{Complex-tangential, finite type domain, interpolation set, pseudoconvex domain}
\subjclass{Primary: 32A38, 32T25; Secondary: 32C25, 32D99}

\begin{abstract} Let $\Omega$ be a bounded, weakly pseudoconvex domain
in $\Cn$, $n\geq 2$, with real-analytic boundary. A real-analytic submanifold 
$\mn \subset \bdy$ is called an analytic interpolation manifold if every real-analytic
function on $\mn$ extends to a function belonging to ${\mathcal{O}}(\overline{\Omega})$. We 
provide sufficient conditions for $\mn$ to be an analytic interpolation manifold. We
give examples showing that neither of these conditions can be relaxed, as well as examples of
analytic interpolation manifolds lying entirely within the set of weakly pseudoconvex points
of $\bdy $.
\end{abstract}

\maketitle

\section{Introduction and Statement of Main Result}

In this paper, we will work with bounded (weakly) pseudoconvex domains $\Omega$
in $\Cn$, $n\geq 2$, with real-analytic boundary. A real-analytic submanifold $\mn$ 
of $\Cn$ contained in $\bdy$ is called an {\bf analytic interpolation manifold} if every 
real-analytic function on $\mn$ extends to some function holomorphic in a neighbourhood of 
$\overline{\Omega}$
(this neighbourhood will, of course, depend on the prescribed function). This definition
is due to Burns and Stout \cite{BS}. Their article proves the following result :

\begin{theorem}[Burns-Stout] Let $\Omega$ be a smoothly bounded strictly pseudoconvex domain
in $\Cn$, $n\geq 2$. A real-analytic submanifold $\mn$ of $\Cn$, $\mn \subset \bdy $, is an analytic 
interpolation manifold if and only if $T_p(\mn ) \subseteq H_p(\bdy ) \ \forall p \in \mn $.
\end{theorem}

In the above result, $H(\bdy )$ is the maximal complex sub-bundle of the tangent bundle $T(\bdy )$.
One could ask whether a real-analytic submanifold $\mn \subset \bdy$, given that $\bdy $ is {\em not}
strictly pseudoconvex along $\mn $, is an analytic interpolation manifold if it is complex
tangential, i.e. if $T_p(\mn ) \subseteq H_p(\bdy ) \ \forall p \in \mn $, and if some
appropriately defined higher Levi-form is strictly positive definite at each point $p \in \mn$.
It will become clear below that complex-tangency is a necessary condition for $\mn $ to be
an analytic interpolation manifold, but the two aforementioned conditions (once precisely defined)
are not sufficient for $\mn $ to be an analytic interpolation manifold. To describe our result,
we need the following definition.

\begin{definition} Let $H^{\Bbb C}(\bdy ) = H^{1,0}(\bdy )\oplus H^{0,1}(\bdy )$ denote the 
complexification of $H(\bdy )$ and $T^{\Bbb C}(\bdy )$ denote the complexification of $T(\bdy )$.
For $p \in \bdy$, let 
${\mathcal{L}}^1_p(\bdy ) = H^{\Bbb C}_p(\bdy )$, and for $j \geq 2$, let
${\mathcal{L}}^j_p(\bdy )$ be the $\Bbb C$-vector space spanned by $H^{\Bbb C}_p(\bdy )$
and all iterated commutators of length $\leq j$ formed by the elements of $H^{\Bbb C}_p(\bdy )$.
A point $p \in \bdy$ is said to be of {\bf Bloom-Graham type $M$} (or simply, of {\bf type $M$}) 
if there exists an $M \in \Bbb N$ such that ${\mathcal{L}}^M_p(\bdy ) = T^{\Bbb C}_p(\bdy )$ 
and ${\mathcal{L}}^j_p(\bdy ) \varsubsetneq T^{\Bbb C}_p(\bdy )$ for $j < M$.
\end{definition}

Our result below shows that if $\bdy $ is of constant type (say $M$) along $\mn$, if $\mn $ is
complex-tangential, and if the $(M-1)$th Levi-form (which is defined below) is positive definite
on a certain subspace of $H_p(\bdy ) \ \forall p \in \mn$, then $\mn $ is an analytic 
interpolation manifold. It is worthwhile noting that the $(M-1)$th Levi form at $p \in \mn$ is 
{\em not} required to be strictly positive definite on all of $H_p(\bdy )$. Furthermore, neither
of the aforementioned conditions can be relaxed. We shall show this through examples in Section 4 
below.
\smallskip

We now define the higher Levi-forms of $\bdy$ that were mentioned above.

\begin{definition} Let $\Omega$ be a smoothly bounded pseudoconvex domain and let 
$p \in \bdy$. Suppose $p$ is
of type $k+1$. Then, we define the \text{\bf $k$th Levi-form of $\bdy$ at $p$},
$ \ \levi{k}_{\bdy }(p \ ;\centerdot) : H^{1,0}_p(\bdy ) \to \Bbb R$ as follows :
\begin{enumerate}
\item There exist holomorphic coordinates $(w_1,...,w_n)$ near $p$ such that 
$\bdy $ is defined, in a neighbourhood of $p$,  by
\begin{equation}
\ws(w) = \sum_{\al + \bet = k+1 \atop 1 \leq \bet < k+1}A^{(p)}_{\alpha\beta}
		{w_*}^{\alpha}{\bar w_*}^{ \ \beta} + {\mathcal{E}}_p(w_*, \mi (w_n))
	 - \er (w_n),
\end{equation}
where we write $w = (w_1,...,w_n) \equiv (w_*,w_n)$ and where ${\mathcal{E}}_p$ is a smooth
function with the property that ${\mathcal{E}}_p(0,0)=0, \ \nabla{\mathcal{E}}_p(0,0)=0$, and
that any term of order $\leq k+1$ is a mixed term involving
non-zero powers of $\mi (w_n), \ w_*$ and ${\bar w}_*$. This is a result from
\cite{BG}. Let $\Fie_p$ be the biholomorphism associated with the above change of coordinate. 
Let $\boo \in H^{1,0}_p(\bdy )$; $d\Fie_p(p)(\boo)$ is an $(n-1)$-tuple 
$d\Fie_p(p)(\boo) = (\zeta_1,...,\zeta_{n-1})$. We define
\[
\levi{k}_{\bdy }(p;\boo) = \sum_{\al + \bet = k+1 \atop  1 \leq \bet < k+1}
				A^{(p)}_{\alpha\beta}{\zeta}^{\alpha}{\bar \zeta}^{\beta}.
\]
\item There is a canonical identification of $H_p(\bdy )$, regarded as a
$\Bbb C$-hyperplane in $\Cn $,  with $H^{1,0}_p(\bdy )$ given by
\begin{equation*}
H_p(\bdy ) \ni (\xi_1,...,\xi_n) \ \rightleftharpoons \ 
		\sum_{l=1}^n \left. \xi_l\Vct{l} \right|_p \in H^{1,0}_p(\bdy ).
\end{equation*}
So, when we say that the Levi-form $\levi{k}_{\bdy }(p \ ;\centerdot)$ acts on 
$(\xi_1,...,\xi_n) \in H_p(\bdy )$, it will mean the action of that
Levi-form on $\sum_{l=1}^n \left. \xi_l\Vct{l} \right|_p \in H^{1,0}_p(\bdy )$.
\end{enumerate}
\end{definition}

We can now state our main result precisely ($\Bbb J$ below is the standard complex
structure map on $\Cn$, and its effect on a vector is equivalent to multiplication
by $i$) :

\begin{theorem} Let $\Omega$ be a bounded, weakly pseudoconvex domain with
real-analytic boundary, and let $\mn$ be a real-analytic, totally real submanifold of 
$\bdy$. Assume that $T_p(\mn ) \subseteq H_p(\bdy )$ for each $p \in \mn$, $\bdy$ is of
constant type $M$ along $\mn$ and that the $(M-1)$th Levi-form of $\bdy$ is positive
definite on the real vector space $ \ {\Bbb J}T_p(\mn ) \subseteq H_p(\bdy )
\ \forall p \in \mn$. Then, $\mn$ is an analytic interpolation manifold.
\end{theorem}

\begin{remark} For $\mn$ to be an analytic interpolation manifold, it is necessary for it to be
totally real, since $\mn$ must not admit any tangential Cauchy-Riemann equations induced by $\bdy$.
This condition on $\mn$ is absent from Theorem 1.1 since it follows from strict pseudoconvexity. 
Additionally, the proof of Theorem 1.4 depends on showing that a complexification of $\mn$ does not
intersect $\overline\Omega$, in the sense of germs, off $\mn$; we use the fact that $\mn$ is totally
real, in our proof, to construct a complexification of $\mn$ that is naturally a complex submanifold
of an open set in $\Cn$.
\end{remark}

\begin{remark} Given $p \in \bdy $, we would like to find a formula for computing 
$\levi{k}_{\bdy }(p \ ;\centerdot)$ directly, without having to find holomorphic
charts near each $p$ in which the defining function has the form (1.1). Also, if $p$ is of 
type $M$, Definition 1.3 only tells us what the $(M-1)$th Levi-form should be, whereas we would
like to be able to compute $\levi{k}_{\bdy }(p \ ;\centerdot)$ for each $k$. Furthermore, we
would like to define $\levi{k}_{\bdy }(p \ ;\centerdot)$ on $H_p(\bdy )$ independently of
the choice of local holomorphic coordinates. These issues are addressed in the next section.
\end{remark}

\section{Preliminary Lemmas}

In this section, we state some general results concerning weakly pseudoconvex domains, which
we shall use in Section 3 to prove Theorem 1.4.

\begin{proposition} Let $\Omega$ and $\mn$ be as in Theorem 1.4. Let $p \in \mn$. There
exist an open neighbourhood $U(p) \subseteq \Cn$ of $p$ and a smooth family 
$\{({\Fie }_q; \omega_q)\}_{q \in U(p) \cap \mn }$ of biholomorphisms 
\[
{\Fie }_q : (\omega_q, q) \to ({\Fie }_q(\omega_q),0)
\]
such that ${\Fie }_q(\omega_q \cap \bdy )$ is defined by
\begin{equation}
\ws_q(w) = \sum_{\al + \bet = M \atop  1 \leq \bet < M}A^{(q)}_{\alpha\beta}
		{w_*}^{\alpha}{\bar w_*}^{ \ \beta} + {\mathcal{E}}_q(w_*, \mi (w_n))
	 - \er (w_n),
\end{equation}
where we write $w = (w_1,...,w_n) \equiv (w_*,w_n)$ and where ${\mathcal{E}}_q$ is a real-analytic
function with the property that ${\mathcal{E}}_q(0,0)=0, \ \nabla{\mathcal{E}}_q(0,0)=0$, and
that any term of order $\leq M$ is a mixed term involving
non-zero powers of $\mi (w_n), \ w_*$ and ${\bar w}_*$. 
\end{proposition}
\begin{proof}
The proof of this statement is standard and originates in \cite{BG}. 
\end{proof}
\smallskip

We now provide a coordinate-free definition for the higher Levi forms, which allows us
to compute $\levi{k}_{\bdy }(p \ ;\centerdot)$ given $p \in \bdy $. To do this, we will
need some preliminary notation. For a multi-index $\alpha = (\alpha_1,...,\alpha_{n-1})
\in {\Bbb N}^{n-1}$ and an integer $1 \leq \mu \leq (n-1)$ define the multi-index $\alpha*\mu$ by
\begin{equation*}
\alpha*\mu = \begin{cases} (\alpha_1,...,\alpha_{\mu-1},\alpha_{\mu}-1,\alpha_{\mu+1},...,\alpha_{n-1}),
			& \text{if $\alpha_{\mu} \geq 1$}\\
		(0,...,0), & \text{otherwise}.
\end{cases}
\end{equation*}
Furthermore, if $\{{\SIl }_j\}_{j=1}^{n-1}$ is any local basis of vector fields for
$H^{1,0}(\bdy )$ near a point $p \in \bdy $, we define ${\SIl }^{\alpha}$ as
\begin{equation*}
{\SIl }^{\alpha} := {\SIl }_{1}^{{\alpha}_1}... \ {\SIl }_{n-1}^{{\alpha}_{n-1}},
\end{equation*}
and we define ${\overline{\SIl }}^{ \ \alpha}$ in an analogous way. Also, we will
use angular brackets $\langle \ , \ \rangle$ to denote contraction between a tangent
vector and a form. 

\begin{definition}[This definition is due to Bloom, \cite{tBl2}]
Let $\Omega$ be a bounded, weakly pseudoconvex domain in 
$\Cn$ with smooth boundary, let $p \in \bdy$, and let $\rho$ be a defining function
for $\Omega$.
\begin{enumerate}
\item An alternate definition of the $k$th Levi-form of $\bdy$ at $p$, 
$ \ \levi{k}_{\bdy }(p \ ;\centerdot) : H^{1,0}_p(\bdy ) \to \Bbb R$ is as follows : 
Let $\{{\SIl }_j\}_{j=1}^{n-1}$ be any local basis of vector fields for $H^{1,0}(\bdy )$ 
near $p$. Then, for any $\boo \in H^{1,0}_p(\bdy )$ 
\begin{equation*}
\levi{k}_{\bdy }(p;\boo ) = \sum_{\al + \bet = k+1 \atop 1 \leq \bet < k+1}
		\dfrac{{\mathfrak{a}}_{\alpha\beta}(p)}{{\alpha}! \ {\beta}!}
				{\zeta}^{\alpha}{\bar \zeta}^{\beta},
\end{equation*}
where $\zeta_j$ are so defined that $\boo = \sum_{j=1}^{n-1}\zeta_j{\SIl }_j|_p$, and where
\begin{equation*}
{\mathfrak{a}}_{\alpha\beta}(p)= -{\SIl }^{\alpha*\mu}\overline{\SIl }^{ \ \beta*\nu}
		\langle[{\SIl }_{\mu}, \overline{\SIl }_{\nu}] \ , \ 
						\overline{\partial}\rho \rangle(p),
\end{equation*}
$\mu$, $\nu$ being so chosen that $\alpha*\mu, \ \beta*\nu \neq 0$ (the coefficient
${\mathfrak{a}}_{\alpha\beta}$ will be independent of the choice of $\mu$, $\nu$ -- this has been
shown in \cite{tBl2}). We note here that $\levi{1}_{\bdy }(p \ ;\centerdot)$ is just the usual 
Levi form of $\bdy $ at $p$ and that if $p$ is a point of type $M$, then
$\levi{k}_{\bdy }(p \ ;\centerdot)=0$ if $k<M-1$.
\smallskip
\item We would like to show that the foregoing definition is the same as Definition 1.3.
This follows from the following result :
\smallskip

\noindent{\bf Theorem.} {\em Let $\Omega$ be as in item (1), $p\in \bdy $, 
and let $p$ be of type $M$. Let $\Fie_p$ be as defined in Proposition 2.1.
The $(M-1)$th Levi-form of ${\Fie}_q({\omega}_q \cap \bdy)$ at the origin,
$\levi{M-1}_{{\Fie}_q({\omega}_q \cap \bdy)}(0 \ ;\centerdot)$ is defined by}
\[
\levi{M-1}_{{\Fie}_q({\omega}_q \cap \bdy)}(0 \ ;\centerdot) : 
{\Bbb C}^{n-1} \ni \zeta \mapsto \sum_{\al + \bet = M \atop 1 \leq \bet < M}A^{(q)}_{\alpha\beta}
		{\zeta}^{\alpha}{\overline{\zeta}}^{\beta}.
\]
\smallskip
The above is a result of Bloom \cite[Theorem 3.3]{tBl2}
\end{enumerate}
\end{definition} 

\begin{remark} Our definition of $\levi{k}_{\bdy }(p \ ;\centerdot)$ (assume that $p$ is a
point of type $M$) differs from that in
\cite{tBl2} by a sign. This is because, in that paper, the normal form for $\bdy $ analogous
to (2.1) above, in local coordinates, is taken to be
\[
\ws(w) = \er (w_n) + \sum_{\al + \bet = M \atop  1 \leq \bet < M}A^{(p)}_{\alpha\beta}
		{w_*}^{\alpha}{\bar w_*}^{ \ \beta} + {\mathcal{E}}_p(w_*, \mi (w_n)).
\]
The reader can check that Definition 2.2(1), when applied to $\levi{1}_{\bdy }(p \ ;\centerdot)$,
gives us the usual Levi-form, which, if $\Omega$ is a pseudoconvex domain, is a positive
semi-definite Hermitian form on $H^{1,0}_p(\bdy )$.
\end{remark}
 
\begin{lemma} Let $\Omega$ be a smoothly bounded pseudoconvex domain. 
Suppose that $0 \in \bdy$ and suppose that, in a neighbourhood of $0$, $\bdy $ is defined by
\[
\rho(w) = h(w_*,\mi (w_n)) - \er (w_n),
\]
where $h$ is a smooth function with $h(0)=0$ and $\nabla h(0)=0$ (and where we write
$w = (w_1,...,w_n) \equiv (w_*,w_n)$).  
Let $D$ be an open neighbourhood of 
$0 \in \Rm$ and let $\gamma=(\gamma_1,...,\gamma_n) : (D,0) \to (\bdy, 0)$ be a smooth
imbedding. Assume that 
\begin{enumerate}
\item $T_{\gamma(x)}(Image(\gamma)) \subseteq H_{\gamma(x)}(\bdy )$ for every $x \in D$;
\item The first Levi-form vanishes on $H^{1,0}_{\gamma(x)}(\bdy )$ for each $\gamma(x)$.
\end{enumerate}
Then, $\gamma_n \equiv 0$.
\end{lemma}
\begin{proof}$\{{\SIl }_j\}_{j=1}^{n-1}$ is a basis of $H^{1,0}(\bdy )$ near the origin, where
${\SIl }_j$ have the form
\[
{\SIl }_j|_{w} = \left. \vect{j}{}\right|_w + \left. F_j{\prr}\vect{n}{}\right|_w.
\]
For any (1,0)-vector field given by ${\Il }|_w = \sum_{j=1}^{n-1}A_j(w){\SIl }_j|_w$, we compute :
\begin{multline}
[{\Il} \ , \ \overline{\Il }] = \ 2i \ \Am\left[\zum\jA\dAl{j}\cvct{\mu}{} + \zum\jA\left\{\dAl{j}\Fl
	+ \Al\dlF{j}\right\}\cvct{n}{}\right. \notag \\ 
	+ \left. \zum\jA\jF\dAl{n}\cvct{\mu}{} + \zum\jA\left\{\Al\jF\dlF{n} +
	\jA\jF\Fl\dAl{n}\right\}\cvct{n}{}\right]. \notag
\end{multline}
Thus, we have 
\begin{equation}
[{\Il} \ , \ \overline{\Il }] = \ \text{\bf V} + 
	2i \ \Am\left[\left(\zum\jA\Al\dlF{j} + \zum\jA\Al\jF\dlF{n}\right)\cvct{n}{}\right],
\end{equation}
where {\bf V} is a section of $H^{\Bbb C}(\bdy )$.

Now consider a (1,0)-vector field, that, restricted to $\text{Image$(\gamma)$}$, is given by
\[
\Il|_{\gamma(x)} = \left. \sum_{j=1}^{n-1}\left\{\sum_{k=1}^m \gamder{j}{k}(x)v_k\right\}\vect{j}{}
		\right|_{\gamma(x)}
	\left. +\left\{\sum_{j=1}^{n-1}\sum_{k=1}^m \gamder{j}{k}(x)F_j[\gamma(x)]v_k \right\}
	\vect{n}{} \right|_{\gamma(x)},
\]
where $(v_1,...,v_m) \in {\Bbb C}^m$.
\[
[\Il \ , \ \overline{\Il}] \ \text{(mod \ $H^{\Bbb C}(\bdy )$)} \equiv 
	\left. \mathfrak{a}(v;x)\cvct{n}{}\right|_{\gamma(x)}-
	\left. \overline{\mathfrak{a}(v;x)}\vect{n}{}\right|_{\gamma(x)} = 0.
\]
The last equality follows from the hypothesis (2) of the lemma. Using (2.2), we get
\[
\mathfrak{a}(v;x)=\sum_{j,\mu=1}^{n-1}\sum_{k,\nu=1}^m 
	\gamder{j}{k}(x)v_k\vect{j}{\overline{F_{\mu}}}[\gamma(x)]
	\overline{\gamder{\mu}{\nu}(x)}\overline{v_{\nu}}
	+\sum_{j,\mu=1}^{n-1}\sum_{k,\nu=1}^m
	\gamder{j}{k}(x)F_j[\gamma(x)]v_k
	\vect{n}{\overline{F_{\mu}}}[\gamma(x)]\overline{\gamder{\mu}{\nu}(x)}
	\overline{v_{\nu}}
\]
In particular, observe that ($\epsilon_{k_0}$ below being the unit vector along the
``$v_{k_0}$-axis'')
\begin{align}
\mathfrak{a}(\epsilon_{k_0};x) &=
	 \zum\vect{j}{\overline{F_{\mu}}}[\gamma(x)]
	\gamder{j}{k_0}(x)\overline{\gamder{\mu}{k_0}(x)}
	+\zum\gamder{j}{k_0}(x)F_j[\gamma(x)]
	\vect{n}{\overline{F_{\mu}}}[\gamma(x)]\overline{\gamder{\mu}{k_0}(x)} \notag \\
	&=\sum_{\mu=1}^{n-1}\fuder{\mu}{k_0}(x)\overline{\gamder{\mu}{k_0}(x)} \notag
\end{align}
since, by complex tangency, $\gamder{n}{k_0}(x) = \sum_{j=1}^{n-1}F_j[\gamma(x)]
\gamder{j}{k_0}(x)$.

Therefore
\[
\mathfrak{a}(v;x) = \ \sum_{k,\nu=1}^m\sum_{\mu=1}^{n-1}\fuder{\mu}{k}(x)
		\overline{\gamder{\mu}{\nu}(x)}v_k\overline{v_{\nu}} \ = \ \langle 
		v|[M_{jk}(x)][D_{jk}(x)]|v\rangle,
\]
where
\begin{align}
M_{jk} &= \fuder{k}{j}(x); \quad j=1,...,m; \ k=1,...,n-1 \notag \\
D_{jk} &= \overline{\gamder{j}{k}(x)}; \quad j=1,...,n-1; \ k=1,...,m \notag \\
\langle v|[A_{jk}]|v \rangle &= \sum_{j,k=1}^m v_jA_{jk}\overline{v_k}. \notag
\end{align}
Consider the sesqui-linear forms
\[
S_x : (u,v) \mapsto  \langle u|[M_{jk}(x)][D_{jk}(x)]|v\rangle; \ x \in D.
\]
Since $\mathfrak{a}(v;x)=0, \ S_x(v,v)=0 \in {\Bbb R}$ for each $v \in {\Bbb C}^m$ and for
each $x \in D$. Thus, $S_x$ are all Hermitian forms that are identically zero. Consequently,
$[M_{jk}(x)][D_{jk}(x)] = 0$ for every $x$. Since $d\gamma(x)$ has maximal
rank for each $x \in D$, we conclude that
\[
\fuder{\mu}{k} \equiv 0; \quad \forall \mu \leq n-1, \ \forall k \leq m.
\]
This implies that, as $F_{\mu}[\gamma(0)] = 0$, 
$F_{\mu}\circ\gamma \equiv 0, \ \forall \mu \leq n-1$. So,
\[
\nabla \gamma_{n}(x)\centerdot v = \sum_{j=1}^{n-1}\sum_{k=1}^m \gamder{j}{k}(x)
	F_j[\gamma(x)]v_k = 0 \quad \forall x \in D, \ \forall v \in {\Bbb R}^m.
\]
This implies that $\nabla {\gamma}_n \equiv 0$, whence ${\gamma}_n \equiv 0$ (since
${\gamma}_n(0)=0$).
\end{proof}
\smallskip

\section{Proof of theorem 1.4}
\begin{proof}[{}]
Without loss of generality, we may assume that $\bdy$ is defined by a global defining function 
$\rho$ that is defined in a neighbourhood $U \supseteq \bdy$. Recall that $\bdy$ is of
type $M$ along $\mn$ and that $T_p(\mn ) \subseteq H_p(\bdy )$ for each $p \in \mn$. 
Pick a $p \in \mn$. There exist a $V(p) \ni p$ open in $\bdy$ and a real-analytic 
imbedding $\gamma : \Rm \supseteq (D,0) \to (\mn \cap V(p),p)$. It must be noted that
$\gamma(\centerdot) \equiv \gamma(\centerdot \ ; p)$, i.e. $\gamma$ depends on $p$, but 
for purposes of notational convenience, we will suppress the dependence on $p$.
It can easily be shown, using standard compactness and homogeneity arguments,
that by our hypothesis on type along $\mn$,
\begin{equation}
\levi{M-1}_{\bdy }(\gamma(x); i(d\gamma(x)\boo)) \geq C|\boo|^{M} \ 
	\forall \boo \in \Rm,
\end{equation}
and for each $p \in \mn$, there exists a neighbourhood $V^*(p) \Subset V(p)$ such that 
(3.1) is true uniformly for all $\gamma(x) \in \text{Image}(\gamma)\cap V^*(p)$ with a 
uniform constant $C \equiv C(p)$.

By Proposition 2.1, there exist an open subset of $\bdy$, $W(p) \Subset V^*(p)$ and a smooth
family of biholomorphisms,
$\{(\Fie_q,{\omega}_q)\}_{q\in W(p)\cap \text{Image}(\gamma)}$, having the effect that for
each $q \in W(p)\cap \text{Image}(\gamma)$, $\Fie_q({\omega}_q\cap\bdy )$ is defined by a
${\ws}_q$ as given in (2.1). Adopting coordinates $w=\Fie_q(z)$, write
\[
\ws_q(w, \bar w) = \peeq(w_*,{\bar w}_*)+ {\mathcal{E}}_q(\mi (w_n),w_*)-\er (w_n),
\]
where $\peeq $ is the polynomial occuring in (2.1).
Let $x_q=\gamma^{-1}(q)$, consider the ball
$B(x_q;\varepsilon_q) \subseteq \gamma^{-1}(\omega_q)$,
and let $\tau_q : x \mapsto (x+x_q)$. Define $\psi_q = \Fie_q\circ\gamma\circ\tau_q$.
Note that $\psi_q : (B(0;\varepsilon_q),0) \to (\Fie_q({\omega}_q\cap\bdy),0)$.
Also $d\psi_q(x) = d\Fie_q(\gamma(x+x_q))\circ d\gamma(x+x_q)$.
\begin{align}
\peeq(i(d\psi_q(0)\boo),-i(\overline{d\psi_q(0)\boo}))
	=& \ \levi{M-1}_{\bdy }(q ; i(d\gamma(x_q)\boo)) \\
	=& \ \levi{M-1}_{\bdy }(\gamma(x_q) ; i(d\gamma(x_q)\boo)) \notag \\
	\geq & \ C|\boo|^M \  \forall \boo \in \Rm, \notag
\end{align}
and from (3.1), we can infer that the above inequality is true uniformly for
$q \in W(p)\cap\text{Image}(\gamma)$. Note that the first equality in (3.2) follows from
Definition 2.2.

Let $\Psi^q$ be the complexification of $\psi_q$ (i.e. $\Psi^q$ is defined, wherever the 
resultant power-series converges, by replacing the real variable $x$ by the complex variable
$\zeta$ in the power-series of $\psi_q$). Since $\{({\Fie }_q;\omega_q)\}_{\rng}$ is a smooth
family, choosing $W(p)$ appropriately, we can find a $\sigma \equiv \sigma(p)$ such that 
$\Psi^q$ are all defined as holomorphic maps on $B(0;\varepsilon_q)+iB(0;\sigma)$ for each
$\rng$. Shrinking $\sigma$ if necessary, we define 
$u_q : B(0;\varepsilon_q)+iB(0;\sigma) \to \Bbb R$ by 
\begin{align}
u_q(\zeta) &= \rho\circ(\Gamma|_{B(x_q;\varepsilon_q)+iB(0;\sigma)})(\zeta+x_q) \notag \\
	&= \ws_q\circ {\Psi}^q(\zeta) \notag
\end{align}
where $\Gamma$ is the complexification of $\gamma$ in an appropriately small 
neighbourhood of $x_q$.

In what follows, we will write $\zeta = \xi + i\eta$, and $\Psi^q \equiv 
(\Psi^q_*,\Psi^q_n)$. By Lemma 2.4, which says that $\Psi^q_n \equiv 0$ 
when $M > 2$, or by the normal form (2.1) in case $M=2$, $u_q$ has the series 
expansion
\begin{align}
u_q(\zeta) =& \sum_{\al + \bet = M \atop 1 \leq \bet < M}A^{(q)}_{\alpha\beta}
		{\Psi^q_*(\zeta)}^{\alpha} \ {\overline{\Psi^q_*(\zeta)}}^{\beta}+
		O(|\zeta|^{(M+1)}) \notag \\
	=& \sum_{\al + \bet = M \atop 1 \leq \bet < M}A^{(q)}_{\alpha\beta}
		\prod_{k=1}^{n-1}{\left(\sum_{l=1}^m\frac{\partial\Psi^q_k}
		{\partial\zeta_l}(0)\zeta_l\right)}^{\alpha_k} \ 
		\prod_{k=1}^{n-1}{\left(\sum_{l=1}^m\overline{\frac{\partial\Psi^q_k}
		{\partial\zeta_l}(0)\zeta_l}\right)}^{\beta_k} + O(|\zeta|^{(M+1)}). \notag
\end{align}
Thus,
\begin{align}
u_q(i\eta) =& \ \peeq(i(d\psi_q(0)\eta),-i(\overline{d\psi_q(0)\eta})) 
		+ O(|\eta|^{(M+1)}) \notag \\
		\geq& \ C|\eta|^M + O(|\eta|^{(M+1)}) && \text{(by (3.2))} \notag
\end{align}
Since the above is true uniformly for all $\rng$ with a uniform constant 
$C \equiv C(p)$, there is a $\delta_p > 0$ such that
\[
u_q(i\eta) > 0; \quad 0<|\eta|< \delta_p, \ \forall\rng .
\]
Another way of saying this is that the complex analytic set 
$\Gamma({\gamma}^{-1}(W(p))+iB(0;\delta_p))$ meets $\overline{\Omega}$ precisely
along $W(p)\cap \mn$.

We can, therefore, find an open neighbourhood $U(p)$ of $p$ in $\Cn$ and a 
complex submanifold $\widetilde{\mn }_p$ of $U(p)$ which is the complexification
of $\mn $ near $p$. $\{U(p)\}_{p\in \mn }$ is an open cover of $\mn$. As $\mn$ is compact,
there exist $p_1,...,p_N \in \mn$ and a tubular neighbourhood $\mathcal{U}$ of $\mn$
such that 
\begin{enumerate}
\item $\mn \subseteq \cup_{k=1}^N U(p_k)$ and $\mathcal{U} \subseteq \cup_{k=1}^N U(p_k)$.
\item $\widetilde{\mn } = \cup_{k=1}^N(\widetilde{\mn }_{p_k} \cap \mathcal{U})$ is a 
complexification of $\mn $, and a complex submanifold of $\mathcal{U}$ such that 
$\widetilde{\mn } \cap \overline{\Omega} = \mn$.
\end{enumerate}
Let $f$ be the real-analytic function prescribed on $\mn$. Shrinking $\mathcal{U}$ if 
necessary, we may assume that $f$ extends to a holomorphic function $\tilde f$
on $\widetilde{\mn }$.

$\overline{\Omega}$ has a basis of Stein neighbourhoods. This follows from
results by Diederich and Fornaess \cite{DF1}, \cite{DF2}, and this is where the assumption
about $\bdy$ being real-analytic gets used. Choose a Stein domain 
$\mathfrak{D} \supset \overline{\Omega}$ such that $\tilde f$ is holomorphic on the 
complex submanifold $(\widetilde{\mn } \cap \mathfrak{D})$ of $\mathfrak{D}$. 
By standard techniques, we can show that $\tilde f$ extends to a 
$F \in \mathcal{O}(\mathfrak{D})$. We remark that this last step reflects a technique
used in \cite{BS} (which follows from theorems A and B of Cartan). Thus, $\mn $ 
is an analytic interpolation manifold.
\end{proof}
\smallskip

\section{Examples}

Before we present our examples, we would like to prove the following proposition.
It can be inferred from \cite[Theorem 1]{BS}, but for the sake of completeness,
we provide a proof.

\begin{proposition} Let $\Omega$ be as in Theorem 1.4. Let 
$\mn$ be a real-analytic submanifold of $\bdy$, let $\widetilde{\mn }$ be
its complexification, and let $p \in \mn $. Suppose there is a curve 
$\gamma \subset \widetilde{\mn }$ passing through $p$ and an $\eps_0 > 0$
such that $[\gamma\cap B(p;\eps)]\cap\overline{\Omega} \neq \emptyset \ \forall\eps
\in (0,\eps_0]$, then $\mn $ is not an analytic interpolation manifold.
\end{proposition}
\begin{proof}
We assume that $\mn$ is an analytic interpolation manifold.
By hypothesis, there exists a real-analytic imbedding $\psi : (S^1,1) \to (\mn ,p)$
onto a simple closed curve $C \subseteq \mn$, so that the following happens :

For $r>0$ sufficiently small, $\psi$ extends to a regular, injective, holomorphic 
map $\Psi$ on $\text{Ann}(0;1-r,1+r)$. There is a small disc $\triangle$, centered at $1$,
such that (defining ${\triangle}^- = (\triangle \cap \{\zeta:|\zeta|<1\})$), without loss
of generality, we have
\begin{gather}
\Psi({\triangle}^-) \cap \overline{\Omega} \neq \emptyset, \notag \\
{\Psi}^{-1}(\Psi({\triangle}^-) \cap \overline{\Omega}) \ \text{contains a curve $L$
tending to $1$}. \notag
\end{gather}

Now choose some $\zeta_0$ in $L$. Define 
\[
g(z) = \frac{1}{(\psi^{-1}(z)-\zeta_0)}
\]
which is real-analytic on $C $, and so extends real-analytically to $\mn $. 
This last conclusion follows from a result by Serre \cite[Secn.19(b)]{jpS}.
By assumption, there exists a $G$ holomorphic in a neighbourhood, call it $D$,
of $\overline{\Omega}$, with $G|_C = g$. We can choose $\triangle$ above to be so small
that $\Psi(\triangle) \subset D$. We can then define
\[
H(\zeta) = G\circ \Psi(\zeta)-\frac{1}{\zeta-{\zeta}_0}.
\]
Clearly $H \in {\mathcal{O}}({\triangle} \setminus \{\zeta_0\})$ and 
$H|_{\triangle \cap S^1} \equiv 0$. The latter implies that $H \equiv 0$. Yet, 
$G\circ \Psi \in {\mathcal{O}}({\triangle})$ whereas $1/(\zeta-{\zeta}_0)$ has a pole 
at $\zeta_0$. This is a contradiction. Our assumption that $\mn$ is an analytic
interpolation manifold must, therefore, be false.
\end{proof}
\smallskip

We can now show that the assumptions on type and positivity in the statement of the
Theorem 1.4 cannot be relaxed. We will do this by constructing real-analytic submanifolds
$\mn$ in $\bdy$ such that $[\widetilde{\mn}]\cap\overline{\Omega} \varsupsetneq \mn$ (here, the
notation $[\widetilde{\mn}]$ denotes the germ of the complexification $\widetilde{\mn}$ along 
$\mn$). In view of Proposition 4.1, $\mn $ would not, therefore, be an analytic interpolation
manifold.
We remark here that if at some $p \in \mn $, $T_p(\mn ) \nsubseteq H_p(\bdy )$, then
$[\widetilde{\mn}]\cap\Omega \neq \emptyset$ near $p$. Complex-tangency is, thus, a {\em necessary}
condition for $\mn$ to be an analytic interpolation manifold.

\begin{example} An example of a domain $\Omega$ and a complex-tangential, totally real submanifold 
$\mn$ where $\bdy$ is of varying type along $\mn$.
\end{example}

Let $\Omega = \{(w_1,w_2) \in {\Bbb C}^2| \ |w_1|^2 +
|w_2|^4 < 1 \}$.

$\Omega$ is a bounded, pseudoconvex domain with real-analytic boundary. $\mn$ will be
a real analytic curve passing through $(-1,0) \in {\Bbb C}^2$, and we will analyze $\mn$
near $(-1,0)$. To simplify calculations, we will work with a biholomorph of $\Omega$. Define
$\omega = \{(z_1,z_2) \in {\Bbb C}^2| \ \er (z_1) > |z_2|^4 \}$. The two domains are
related by a biholomorphism
\[
\varPhi(z_1,z_2) = (w_1,w_2) = 
	\left(\dfrac{z_1-1}{z_1+1}, \dfrac{\sqrt2 z_2}{\sqrt{z_1+1}}\right),
\]
where $\varPhi : (\omega,\BD ) \to (\Omega,(\bdy \setminus \{(1,0)\}))$ ({\bf Note :} Since
$\er (z_1) > 0$ when $(z_1,z_2) \in \overline{\omega}$, we can choose an appropriate
analytic branch of $z_1 \mapsto \sqrt{z_1+1}$.). Consider the 
real-analytic, complex-tangential curve $\gamma : \Bbb R \to \BD$ (and write
$\MN = \text{Image}(\gamma)$)
\[
\gamma(t) = (t^4,t).
\]
Define $\mn = \overline{\varPhi(\MN )}$. It is easy to check that $\mn \subseteq \bdy $ is a
real-analytic, complex-tangential submanifold. Notice that at $t=0$, $\text{Image}(\gamma)$ 
passes through a point of type $4$, whereas it is of type $2$ elsewhere. Consequently, the 
Bloom-Graham type of $\bdy $ varies along $\mn$.

Let $\Gamma$ denote the complexification of $\gamma$. Writing $\zeta = Re^{i\theta}$
\begin{align}
\er [{\Gamma}_1(\zeta)] &= \ R^4\cos4\theta, \notag \\
|{\Gamma}_2(\zeta)|^4 &= \ R^4. \notag
\end{align}
Observe that $\er [{\Gamma}_1(Re^{i\theta})]=|{\Gamma}_2(Re^{i\theta})|^4$
for $\theta = 0, \ \pi/2, \ \pi$ and $3\pi/2$. Thus
$\boldsymbol{[}\widetilde{\MN }\boldsymbol{]} \cap \overline{\omega} \varsupsetneq \MN$. 
Since the notion of type is invariant under biholomorphic transformations, we have found
a real-analytic, complex-tangential $\mn \subseteq \bdy$ such that $\bdy$ is of varying
type along $\mn$ and such that 
$\boldsymbol{[}\widetilde{\mn }\boldsymbol{]} \cap \overline{\Omega} \varsupsetneq \mn$.
\medskip

We comment on the notation used in the next three examples. In all of the equations
involving Levi-forms, we will use the identification between $H^{1,0}_p(\bdy )$ and
$H_p(\bdy )$ that was introduced in Definition 1.3(2).

\begin{example} An example of a domain $\Omega$ and a complex-tangential, totally real submanifold 
$\mn$ where $\bdy$ is of constant type along $\mn$, but where the positivity condition fails.
\end{example}

Let $\Omega = \{(w_1,w_2,w_3,w_4) \in {\Bbb C}^4| \ |w_1|^2 + \sum_{k=2}^4 |w_k|^4 < 1 \}$.

$\Omega$ is a bounded, pseudoconvex domain with real-analytic boundary. As in Example 4.2, we
will work with a biholomorph of $\Omega$. Define 
$\omega = \{(z_1,z_2,z_3,z_4) \in {\Bbb C}^4| \ \er (z_1) > \sum_{k=2}^4 |z_k|^4 \}$.
The two domains are related by a biholomorphism
\[
\varPhi(z_1,z_2,z_3,z_4) = (w_1,w_2,w_3,w_4) = 
	\left(\dfrac{z_1-1}{z_1+1}, \dfrac{\sqrt2 z_2}{\sqrt{z_1+1}}, \dfrac{\sqrt2 z_3}
			{\sqrt{z_1+1}},\dfrac{\sqrt2 z_4}{\sqrt{z_1+1}} \right),
\]
where $\varPhi : (\omega,\BD ) \to (\Omega,(\bdy \setminus \{(1,0,0,0)\}))$. Consider the 
real-analytic, complex-tangential curve $\gamma : [-\pi,\pi] \to \BD$ (and write 
$\text{Image}(\gamma) = \MN$)
\[
\gamma(\theta) = ({\sin}^4\theta+{\cos}^4\theta+1,\ 1,\ \sin\theta,\ \cos\theta).
\]
Let $\mn = \varPhi(\MN )$. $\mn \subseteq \bdy $ is clearly a
real-analytic, complex-tangential submanifold.

Notice that $H(\BD )$ is spanned at every point by the vector fields
\begin{gather}
{\Il }_1|_{(z_1,z_2,z_3,z_4)} = (4\overline{z_2}|z_2|^2,1,0,0), \qquad\qquad 
			{\Il }_2|_{(z_1,z_2,z_3,z_4)} = (4\overline{z_3}|z_3|^2,0,1,0), \notag \\
{\Il }_3|_{(z_1,z_2,z_3,z_4)} = (4\overline{z_4}|z_4|^2,0,0,1). \notag
\end{gather}
Observe further that the complex Hessian of $\rho$ (where $\rho$ is the defining function
of $\BD $) is given by
\[
({\mathfrak{H}}_{\Bbb C}\rho)(z_1,z_2,z_3,z_4) = \begin{pmatrix}
					0 & 0 & 0 & 0\\
					0 & 4|z_2|^2 & 0 & 0\\
					0 & 0 & 4|z_3|^2 & 0\\
					0 & 0 & 0 & 4|z_4|^2
					  \end{pmatrix}.
\]
Let $\levi{1}_{\BD }(p \ ;\centerdot)$ denote the Levi-form at $p \in \BD$. Notice that
$\levi{1}_{\BD }(\centerdot \ ; {\Il }_1) \neq 0$ along $\MN$. Thus $\BD$ is of constant
Bloom-Graham type $2$ along $\MN$. This implies that $\bdy$ is of constant Bloom-Graham type $2$
along $\mn$.
But
\[
\levi{1}_{\BD }(\gamma(\theta); \ i{\gamma}^{\prime}(\theta)) = 8{\cos}^2\theta \ {\sin}^2\theta
\]
which vanishes at $\gamma(0)$, whence the positivity condition along $\MN $ fails.

Let $\Gamma$ denote the complexification of $\gamma$. Writing $\zeta = \xi+i\eta$, we observe
\begin{align}
\er [{\Gamma}_1(i\eta)] &= \ {\cosh}^4\eta + {\sinh}^4\eta + 1, \notag \\
|{\Gamma}_2(i\eta)|^4 + |{\Gamma}_3(i\eta)|^4 + |{\Gamma}_4(i\eta)|^4 &= \ 
{\cosh}^4\eta + {\sinh}^4\eta + 1. \notag
\end{align}
Thus $\boldsymbol{[}\widetilde{\MN }\boldsymbol{]} \cap
\overline{\omega} \varsupsetneq \MN$, whence, arguing exactly as in Example 4.2 above,
we have a real-analytic, complex-tangential $\mn \subseteq \bdy$ such that $\bdy$ is of
constant type along $\mn$, that the positivity condition fails and such that
$\boldsymbol{[}\widetilde{\mn }\boldsymbol{]} \cap \overline{\Omega} \varsupsetneq \mn$.
\medskip

We would now like to give examples of analytic interpolation manifolds. In Example 4.4,
our manifold $\mn $ passes through weakly pseudoconvex points. In Example 4.5, $\mn $
runs through points of Bloom-Graham type 4, although $\levi{3}_{\bdy }(p \ ;\centerdot)$
is not strictly positive definite on $H^{1,0}_p(\bdy )$ for any $p \in \mn$.

\begin{example} An example of a weakly pseudoconvex domain and an analytic
interpolation manifold.
\end{example}

Let $\Omega$ be exactly as in Example 4.3. All the notation used below will have the same
meanings as in Example 4.3. As in that example, we will work with
\[
\omega = \{(z_1,z_2,z_3,z_4) \in {\Bbb C}^4| \ \er (z_1) > \sum_{k=2}^4 |z_k|^4 \},
\]
a biholomorph of $\Omega$. Consider the real-analytic, complex-tangential curve
$\gamma : [-\pi,\pi] \to \BD $ (and write $\text{Image}(\gamma) = \MN$)
\[
\gamma(\theta)= ((2+\cos\theta)^4+(2+\sin\theta)^4,\ 0,\ 2+\sin\theta,\ 2+\cos\theta).
\]
As before, define $\mn = \varPhi(\MN )$, which is a complex-tangential, real-analytic
submanifold of $\bdy$.

As in Example 4.3 let $\levi{1}_{\BD }(p \ ;\centerdot)$ denote the Levi-form at $p \in \BD$.
\[
\levi{1}_{\BD }(\gamma(\theta); \ {\Il }_3|_{\gamma(\theta)})=4(2+\cos\theta)^2 \neq 0,
\]
whence $\BD$ is of constant type $2$ along $\MN$. Consequently, $\bdy$ is of type $2$ along
$\mn$, although {\em every} point on $\mn$ is a weakly pseudoconvex point. To see this, 
observe that $\levi{1}_{\BD }(\centerdot \ ; {\Il }_1) = 0$ along $\MN$. But
\[
\levi{1}_{\BD }(\gamma(\theta); \ i{\gamma}^{\prime}(\theta)) = 
4{\cos}^2\theta(2+\sin\theta)^2+4{\sin}^2\theta(2+\cos\theta)^2 > 0.
\]
Since biholomorphisms preserve the positivity of the Levi-form and since $\varPhi$ is biholomorphic
in a neighbourhood of $\MN $, our positivity condition is
preserved for $\mn$. So $\mn \subseteq \bdy$ is a complex-tangential, real-analytic submanifold
of $\bdy$ that satisfies our positivity condition. By Theorem 1.4, therefore, $\mn$ is an 
analytic interpolation manifold.
\medskip

\begin{example} Another example of a weakly pseudoconvex domain and an analytic 
interpolation manifold $\mn$. In this example, each point of $p \in \mn $ is a point of type 
4 and, in fact, $\levi{3}_{\bdy }(p \ ;\centerdot)$ is negative in certain directions
in $H^{1,0}(\bdy )$.
\end{example}

Let $\Omega = \{(z,w) \in {\Bbb C}^2 | \ |w+e^{i\lgz}|^2+C[\lgz ]^4 < 1 \}$.

This example is taken from \cite{aN}. For an appropriate $C > 0$, $\Omega$ is a 
pseudoconvex domain. We define
\[
\mn = \{(z,w) \in \bdy | \ |z|=1, w = 0 \}.
\]
$\mn $ has a real-analytic parametrization $\gamma : [-\pi,\pi] \to \bdy $ given by
$\gamma(\theta) = (e^{i\theta},0)$. As in \cite{aN}, we can show that each point of
$\mn $ is of type $> 2$.

In what follows, we will write $A(z) = e^{i\lgz}$, $B(z) = e^{-i\lgz}$. 
$H^{1,0}(\bdy )$ is spanned by the vector field
\[
{\Il }|_{(z,w)} = \ z[{\bar w}+B(z)] \ \Vct{}-
			[i{\bar w}A(z)-iwB(z)+4C\{\lgz \}^3] \ \vect{}{}.
\]
We  can now see that $\gamma^{\prime}(\theta) \in H_{(e^{i\theta},0)}(\bdy)$. In fact,
by the identification introduced in Definition 1.3(2)
\[
i\gamma^{\prime}(\theta) \ \rightleftharpoons \ -{\Il }|_{(e^{i\theta},0)}.
\]
We get the equation
\begin{multline}
\levi{3}_{\bdy }((e^{i\theta},0); \ i{\gamma}^{\prime}(\theta))\\
	= \ -\left\{(-1)^3(-1)\dfrac{{\Il }^2}{3!}+(-1)^2(-1)^2
			\dfrac{{\Il }\overline{\Il }}{2! \ 2!}+
			(-1)(-1)^3\dfrac{\overline{\Il }^2}{3!}\right\}
	\langle[{\Il },\overline{\Il }] \ , \overline{\partial}\rho \rangle
	(e^{i\theta},0).
\end{multline}
We compute to find that
\begin{equation}
-\langle[{\Il },\overline{\Il }] \ , \overline{\partial}\rho \rangle(z,w)
 = \ 12C[\lgz ]^2 + O([\lgz ]^3, \ wB(z), \ {\bar w}A(z), \ |w|^2, \ |w|[\lgz ]^2).
\end{equation}
It can easily be shown from (4.1) and (4.2) that
\[
\levi{3}_{\bdy }((e^{i\theta},0); \ i{\gamma}^{\prime}(\theta)) = 14C > 0.
\]

Similarly, we can show that 
$\levi{3}_{\bdy }((e^{i\theta},0); \ i\Il ) = -2C < 0$. 
So, $\levi{3}_{\bdy }((e^{i\theta},0) \ ; \centerdot)$ is actually negative in certain
directions in $H^{1,0}_{(e^{i\theta},0)}(\bdy )$. Yet, by Theorem 1.4, $\mn $ is an analytic
interpolation manifold (in this case, we can also check very easily that
$\boldsymbol{[}\widetilde{\mn }\boldsymbol{]} \cap \overline{\Omega}=\mn$).
\bigskip

\noindent{\bf Acknowledgement.} The author wishes to thank Alexander Nagel for his 
encouragement and for the many useful discussions during the course of this work.


\smallskip

\begin{thebibliography}{7}

\bibitem{tBl2}
T. Bloom, {\em Remarks on type conditions for real hypersurfaces in 
${\Bbb C}^n$}, Several Complex Variables - Proceedings of International 
Conferences, Cortona, Italy 1976-1977, Scuola Norm. Sup. Pisa, 1978, pp. 14-24.

\bibitem{BG}
T. Bloom and I. Graham, {\em On type conditions for generic real 
submanifolds in ${\Bbb C}^n$}, Invent. Math. {\bf 40} (1984) 217-243.

\bibitem{BS}
D. Burns and E.L. Stout, {\em Extending functions from submanifolds of
the boundary}, Duke Math. J. {\bf 43} (1976), 391-404.

\bibitem{DF1}
K. Diederich and J.E. Fornaess, {\em Pseudoconvex domains : existence of 
Stein neighbourhoods}, Duke Math. J. {\bf 44} (1977), 641-662.

\bibitem{DF2}
K. Diederich and J.E. Fornaess {\em Complex submanifolds in real-analytic
pseudoconvex hypersurfaces}, Proc. Nat. Acad. Sci. {\bf 74} (1977), 3126-3127.

\bibitem{aN}
A. Noell {\em Properties of peak sets in weakly pseudoconvex domains in ${\Bbb C}^2$},
Math. Z. {\bf 186} (1984), 99-116.

\bibitem{jpS}
J.-P. Serre, {\em Applications de la th\'eorie g\'en\'erale \`a divers
probl\`emes globaux}, S\'eminaire H. Cartan, 1951-52.
\end{thebibliography}
\enddocument